\documentclass[12pt]{article}

\usepackage{amsfonts}
\usepackage{xr}
\usepackage{graphicx}
\usepackage{amsmath}

\def\2{\mathcal{C}^2(\mathbb{R}^N)}
\def\A{\mathcal{A}}

\def\R{\mathbb{R}}

\def\disp{\displaystyle}
\def\tilde{\widetilde}
\def\epsilon{\varepsilon}
\def\trait (#1) (#2) (#3){\vrule width #1pt height #2pt depth #3pt}
\def\fin{\hfill\trait (0.1) (5) (0) \trait (5) (0.1) (0) \kern-5pt \trait (5) (5) (-4.9) \trait (0.1) (5) (0)}

\newtheorem{thm}{\bf Theorem}[section]

\newtheorem{prop}[thm]{\bf Proposition}

\newtheorem{defi}[thm]{\bf Definition}

\title{The effect of the Schwarz rearrangement on the periodic principal eigenvalue of a nonsymmetric operator}
\author{Gr\'egoire Nadin\thanks{
D\'epartement de Math\'ematiques et Applications,
\'Ecole Normale Sup\'erieure, CNRS UMR8553 ,
            45 rue d'Ulm, F~75230 Paris cedex 05, France ; e-mail : nadin@dma.ens.fr}}

\begin{document}
\date{}
\maketitle

\begin{abstract}
This paper is concerned with the periodic principal eigenvalue $k_{\lambda}(\mu)$ associated with the operator
\begin{equation}-\frac{d^2}{dx^2} -2\lambda \frac{d}{dx} -\mu(x)-\lambda^2,\end{equation}
where $\lambda\in\R$ and $\mu$ is continuous and periodic in $x\in\R$. Our main result is that $k_\lambda(\mu^*)\leq k_\lambda(\mu)$, 
where $\mu^*$ is the Schwarz rearrangement of the function $\mu$. From a population dynamics point of view, using reaction-diffusion modeling, 
this result means that the fragmentation of the habitat of an invading population slows down the invasion. We prove that this property does not hold in 
higher dimension, if $\mu^*$ is the Steiner symmetrization of $\mu$. For heterogeneous diffusion and advection, we prove that increasing the period of the 
coefficients decreases $k_\lambda$ and we compute the limit of $k_\lambda$ when the period of the coefficients goes to $0$. Lastly, we prove that, in dimension $1$, rearranging the diffusion 
term decreases $k_\lambda$. 
These results rely on some new formula for the periodic principal eigenvalue of a nonsymmetric operator.
\end{abstract}


\bigskip

\bigskip

\noindent {\bf Abbreviated title:} Schwarz rearrangement of a nonsymmetric operator

\section{Introduction and main results}

\subsection{General framework and hypotheses}

This paper is mainly concerned with the periodic principal eigenvalue of a nonsymmetric operator $L_\lambda$. This operator is defined for all $\phi\in\mathcal{C}^2(\R)$ by
\begin{equation} \label{defLlambda}-L_\lambda\phi=-\frac{d^2}{dx^2}\phi -2\lambda \frac{d}{dx}\phi -(\mu(x)+\lambda^2)\phi,\end{equation}
where $\lambda\in\R$. 
The zero order term $\mu$ is supposed to be continuous and periodic in $x$, namely, there exist some positive constants $L$ such that for all $x\in\R$:
\begin{equation} 
\mu(x+L)=\mu(x).
\end{equation} 
The period $L$ will be fixed in the sequel and when
a function is said to be periodic, this periodicity
will always refer to this given period (except in section \ref{sectionperiod}). We thus denote $\mathcal{C}^0_{per}(\R)$ (resp.  $\mathcal{C}^1_{per}(\R)$) the set of the periodic functions that are continuous (resp. of class $\mathcal{C}^1$). 

The \textit{periodic principal eigenvalue} $k_{\lambda}=k_{\lambda}(\mu)$ of the operator $L_\lambda$ is uniquely defined by the existence of a \textit{periodic principal eigenfunction} $\phi\in\mathcal{C}^2(\R)$ (which is unique up to multiplication by a positive constant) that satisfies:
\begin{equation} \label{eigen0} \left\{ \begin{array}{l}
-L_{\lambda} \phi = k_\lambda \phi, \\
\phi > 0, \ \phi \hbox{ is periodic}.\\
\end{array} \right.
\end{equation}

The existence and uniqueness of these eigenelements have been proved in \cite{Pinsky} for example. We are concerned here with the dependence relation between $\mu$ and $k_\lambda(\mu)$.

\subsection{The link with reaction-diffusion equations}

This periodic principal eigenvalue is involved in
the qualitative properties of the solutions of the heterogeneous reaction-diffusion equation
\begin{equation} \label{eqprinc} 
 \partial_t u -\Delta u=f(x,u) \hbox{ for all } \ (t,x)\in\R\times\R^N.
\end{equation}

This equation arises in ecology, genetics, chemistry or combustion models. We will focus here on the ecology modeling, which has been first discussed by J. G. Skellam \cite{Skellam} in 1951. In this framework, $u$ represents a population density and $f$ is a growth rate. The first rigorous mathematical results for the homogeneous equation in a bounded domain with Dirichlet boundary conditions, are due to D. Ludwig, D. Aronson and H. Weinberger \cite{LudwigAronsonWeinberger}. They proved that, under certain hypotheses on $f$, the sign of the principal eigenvalue associated with the linearization of the equation in a neighborhood of $0$ determines the existence and the attractivity of a non-null steady state. This property has been generalized to various heterogeneous equations in bounded domains \cite{CantrellCosner1, CantrellCosner2, CantrellCosnermathsbio, MurraySperb}.

\bigskip

We are interested in the present paper in unbounded domains. The reason for this choice will be given later. The case of a general unbounded domains is difficult to investigate since the principal eigenvalue of a linear operator does not exist in general. This is why we assume that the growth rate is periodic in $x$. In general unbounded domains, one needs to introduce the notion of generalized principal eigenvalues, defined by H. Berestycki, F. Hamel and L. Rossi in \cite{Rossi}. We will investigate the extension of properties proved in this paper in a forthcoming work. This extension is far from being obvious since, in unbounded media, several notions of generalized principal eigenvalues can be used and it is not clear if these notions are equivalent or not. 

Equation (\ref{eqprinc}) with periodic boundary conditions has first been studied by H. Berestycki, F. Hamel and L. Roques \cite{Base1}. They made the following hypotheses on $f$:
\begin{equation} \label{hypf} \left\{ \begin{array}{l} 
f(x,0)=0 \hbox{ for all} \ x\in\mathbb{R}^{N},\\
\forall x \in \mathbb{R}^{N}, \ s \rightarrow
f(x,s)/s \ \hbox{decreases on} \ \mathbb{R}^{+*},\\
\exists \ M>0 \ | \ \forall x \in \mathbb{R}^{N}, \forall s \geq M, \ f(x,s) \leq 0.\\ 
\end{array} \right. \end{equation}
These hypotheses make sense from an ecological point of view. The function $s\mapsto f(x,s)/s$ is a growth rate per capita and it decreases because of the intraspecific competition for ressources. As the carrying capacity of the media is limited, the growth rate $f(x,u)$ will be negative when the population density is too large.

H. Berestycki, F. Hamel and L. Roques proved that equation (\ref{eqprinc}) admits a positive periodic solution if and only if $k_0(\mu)<0$, where $\mu(x)=f_u'(x,0)$. Moreover, this solution is unique when it exists and attracts all the solutions of the associated Cauchy problem with nonnegative and nonnull initial datum. Thus, the negativity of the periodic principal eigenvalue $k_0(\mu)$ is a necessary and sufficient condition for the survival of the population. It is then natural to study the dependence relation between the coefficient $\mu$ and $k_0(\mu)$. 
From a biological point of view, it is interesting to study the influence of the heterogeneity of the medium on the survival of the population. Various quantifications of this ``heterogeneity'' are given in \cite{Base1}, which all yield that this heterogeneity is better for the survival of the population. 

\bigskip

But, when the population survives, these results do not quantify ``how'' it survives. If one considers an invading population, it is interesting to know if this population can survive, but also the speed of its invasion. In order to compute this speed, the biologists try to find the level lines of the population density (see \cite{Skellam} for example). Some empirical observations yield that the square root of the level sets areas is a linear function of the time. The associated linear growth rate can thus be viewed as an invasion speed.

Such an observation is asymptotic in $t$ and thus one needs to consider unbounded domains in order to modelize this phenomenon. This is why we are interested in periodic boundary conditions in the present paper. 

A mathematical modeling for this empirical observation has been given by D. Aronson and H. Weinberger \cite{Aronson}, when $f$ does not depend on $x$, and by M. Freidlin and J. Gartner \cite{Gartner}, when $N=1$, $f$ is positive between $0$ and $1$, $f(x,1)=0$ for all $x$ and $f(x,s)\leq f_u'(x,0)s$ for all $(x,s)$. They proved that there exists $w^*_e>0$ for $e=\pm 1$ such that for all solution $u$ of (\ref{eqprinc}) with compactly supported, non-null and nonnegative initial datum, one has:
\begin{equation} \left\{ \begin{array}{rcccl}
u(t,x+wte)\rightarrow 0 &\hbox{ as }& t\rightarrow+\infty &\hbox{ for all }& x\in\R^N, \ w>w^*_e,\\
u(t,x+wte)\rightarrow 1 &\hbox{ as }& t\rightarrow+\infty &\hbox{ for all }& x\in\R^N, \ 0\leq w<w^*_e.\\
\end{array} \right. \end{equation}
This property means that the level set $\Gamma_t=\{x, \ u(t,x)=\frac{1}{2}\}$ ``looks like'' $\{-w^*_{-e}t,w^*_e t\}$ as $t\rightarrow+\infty$. 
The quantity $w^*_e$ is thus called the \textit{spreading speed} in direction $e$. 
It can be characterized with the help of the family of the periodic principal eigenvalues $(k_\lambda)_{\lambda>0}$ as 
$w^*_e=\min_{\lambda>0} -k_{\lambda e}/\lambda.$

This property has been extended to a general multidimensional framework by H. Weinberger in \cite{Weinberger}, using discrete dynamical systems tools, and H. Berestycki, F. Hamel and G. Nadin in \cite{spreadinggen}, using PDE's tools. In multidimensional media, when $f(x,s)\leq f_u'(x,0)s$ for all $(x,s)$, the spreading speed is characterized for all $e\in\mathbb{S}^{N-1}$ as
\begin{equation} \label{formulew^*} w^*_e=\min_{\lambda\cdot e>0} \frac{-k_\lambda}{\lambda\cdot e},\end{equation}
where $\lambda\in\R^N$ and $k_\lambda$ is the periodic principal eigenvalue associated with 
the operator $L_{\lambda}$, defined by
\begin{equation}\label{LlambdamultiD} -L_\lambda\phi=-\Delta\phi -2\lambda \cdot\nabla\phi - (|\lambda |^2 +\mu(x))\phi.\end{equation}

The notion of \textit{pulsating traveling fronts} has been given in a parallel way by N. Shigesada, K. Kawasaki and E. Teramoto in \cite{Shigesada1} and by J. Xin in \cite{Xin}. If $f$ is positive between $0$ and $1$ and $f(x,1)=0$ for all $x$, a pulsating traveling front of equation (\ref{eqprinc}) of speed $c$ in direction $e$ that connects $0$ to $1$ is a solution $u$ that can be written $u(t,x)=\phi(x\cdot e-ct,x)$, where $\phi=\phi(z,x)$ is periodic in $x$, $\phi(-\infty,x)=1$ and $\phi(+\infty,x)=0$. The existence of pulsating traveling fronts for equation (\ref{eqprinc}) has been proved in various contexts \cite{Frontexcitable, Base2, twperiodic, Xin}. If $f(x,s)\leq f_u'(x,0)s$ for all $(x,s)$, then there exists a pulsating front of speed $c$ in direction $e$ if and only if $c\geq c^*_e$, where
\begin{equation} \label{formulec^*} c^*_e=\min_{\lambda >0} \frac{-k_{\lambda e}}{\lambda}. \end{equation}
In dimension $1$, $w^*_e=c^*_e$ if $e=\pm 1$ and we expect that for all nonnegative, non-null initial datum, the solution of the associated Cauchy problem converges to a pulsating traveling front along its level lines. In multidimensional media, $w^*_e$ can be computed with the help of the $(c^*_\xi)_{\xi\in\mathbb{S}^{N-1}, e\cdot \xi>0}$. 

\bigskip

Thus the speed $w^*_e$ quantifies the speed of a biological invasion and it is relevant to investigate the influence of the heterogeneity of the media on this speed. For example, in \cite{Base2}, H. Berestycki, F. Hamel and L. Roques have proved that $c^*_e(\overline{\mu})\leq c^*_e(\mu)$ for all $\mu, e$, where $\overline{\mu}$ is the average of $\mu$ in the periodicity cell, or that $B\mapsto c^*_e(B\mu)$ is nondecreasing if $\int_C \mu\geq 0$. This proves that the more the amplitude of a given heterogeneity is large, the more the invasion speed is high. 

It is more difficult to determine the influence of the shape of the heterogeneity on the speed. Consider for example what the biologists call a ``patchy environment'', that is, a growth rate $\mu$ which can be written
\begin{equation}\label{patchy} \mu(x)= \left\{\begin{array}{rcl}
\mu^+ \ &\hbox{ if }& x\in\Omega,\\
\mu^- \ &\hbox{ if }& x\notin\Omega,\\
\end{array} \right. \end{equation}
with $\mu^+>\mu^-$ and $\Omega$ is a measurable set. In this case, the set $\Omega$ is called the \textit{habitat} since it is the place where the growth rate is high. It is relevant from a biological point of view to try to determine how the fragmentation of the habitat decreases the invasion speed. 

But from a mathematical point of view, it is difficult to quantify the fragmentation of a set. It seems reasonable to say that the more the period of the environment is large, the more the environement is fragmented. 
In \cite{Shigesada2, Shigesada1}, N. Shigesada, K. Kawazaki and co-authors numerically observed that if $\mu_L(x)=\mu(x/L)$, with $\mu$ as in (\ref{patchy}), then $L\mapsto c^*_e(\mu_L)$ is nondecreasing for all $e$. But then we only rescale the set and, somehow, the ``shape'' of the set does not vary.

In dimension $1$, in \cite{Base1} H. Berestycki, F. Hamel and L. Roques noticed that the less fragmented measurable sets $\Omega$ of a fixed measure $m$ are the interval $\Omega=(0,m)$ and all its translations. They also generalized this kind of construction to other growth rates $\mu$ using the notion of the Schwarz rearrangements.


\subsection{Rearrangements and eigenvalues optimization}

The classical rearrangement for periodic functions of the real line is the Schwarz rearrangement. For more details and properties about
this notion, we refer to \cite{Kawohl} for example.

\begin{defi}\cite{Kawohl}
 Assume that $\mu\in \mathcal{C}^0(\R)$ is a L-periodic function. Then there exists a unique L-periodic function $\mu^{*}\in \mathcal{C}^0(\R)$ such that:

i) $\mu^{*}$ is symmetric with respect to $L/2$,

ii) $\mu^{*}$ is nondecreasing on $(0,L/2)$,

iii) $\mu^{*}$ has the same distribution function as $\mu$, that is for all $\alpha\in\R$:
\[meas \{ x, \mu(x)>\alpha \} = meas \{ x, \mu^{*}(x)>\alpha \}.\]

The function $\mu^{*}$ is called the \textit{Schwarz periodic rearrangement} of the function $\mu$.
\end{defi}

For example, if $\mu=1_{\Omega}$ in $(0,L)$, where $\Omega$ is a measurable set of measure $m$, then $\mu^*$ is the characteristic function of the set $(\frac{L-m}{2},\frac{L+m}{2})$. This example illustrates how the Schwarz rearrangement of a growth rate $\mu$ is, somehow, the growth rate associated with the less fragmented environment, when the distribution function of $\mu$ is fixed. 

\bigskip

In \cite{Base1}, H. Berestycki, F. Hamel and L. Roques proved that $k_0(\mu^*)\leq k_0(\mu)$ for all $\mu$. This means that, somehow, the fragmentation of the habitat disadvantages the survival of the population since the sign of $k_0(\mu)$ determines the existence and the attractivity of a non-null equilibrium state for equation (\ref{eqprinc}) under some hypotheses. We can wonder if it also increases the invasion speed or, equivalently, if $k_{\lambda}(\mu^*)\leq k_{\lambda}(\mu)$ for all $\lambda$. This was left as an open question by H. Berestycki, F. Hamel and L. Roques in \cite{Base2}. The main result of the present paper answers this question:

\begin{thm} \label{influencesteiner}
If $N=1$, then for all $\mu\in\mathcal{C}^0_{per}(\R)$, one has:
\[k_\lambda(\mu^*)\leq k_\lambda(\mu).\]
\end{thm}

As a corollary, we immediatly get $c^*_e(\mu^*)\geq c^*_e(\mu)$ for $e=\pm 1$.

The proof of $k_0(\mu^*)\leq k_0(\mu)$ in \cite{Base1} strongly relies on the Rayleigh characterization of the periodic principal eigenvalue of a self-adjoint operator:
\begin{equation} \label{RayleighSteiner} k_0(\mu)=\min_{\alpha\in \mathcal{C}^1_{per}(\R)}\frac{\int_0^L \alpha'^2-\int_0^L\mu\alpha^2}{\int_0^L \alpha^2}.\end{equation}

The main difficulty with the proof of Theorem \ref{influencesteiner} is that $k_{\lambda}(\mu)$ is the periodic principal eigenvalue of a nonsymmetric operator $L_\lambda$. Thus the Rayleigh characterization (\ref{RayleighSteiner}) does not hold anymore and one has to find a new method. 

\bigskip

In the case of Dirichlet, Neumann or Robin boudary conditions, there is a wide litterature on the link between rearrangements and eigenvalues (see \cite{Henrot, HenrotPierre} for example). But very few papers consider nonsymmetric operators and the Rayleigh characterization is always the main tool that is used. Only F. Hamel, N. Nadirashvili and E. Russ \cite{FaberKrahndrift, isoperimetric, HamelNadirashviliruss} gave optimization results for the principal eigenvalue of the operator $-\Delta-q(x)\cdot \nabla-\mu$ with Dirichlet boundary conditions. A. Alvino, G. Trombetti and P.-L. Lions  also investigated similar problems in \cite{Lionsdrift, Lions}. But these works enable some modification of the first order term $q$. In our case, $q$ is a constant $2\lambda$ and we want it to stay in this form.

Lastly, when the diffusion term or advection term are heterogeneous, hardly no paper have previously investigated the influence of the shape of this heterogeneity. We prove at the end of this paper that increasing the period of the coefficients decreases the periodic principal eigenvalue in any dimension and that rearranging the diffusion term also decreases this eigenvalue in dimension $1$. 


\subsection{Counterexample in higher dimension}

Consider now a function $\mu$ periodic on the set $\R^{N}$. That is, there exist some positive constants $L_{1},...,L_{N}$ such that for all $x,i$:
\begin{equation} 
\mu(x+L_{i}\epsilon_i)=\mu(x),
\end{equation}
where $(\epsilon_i)_{i\in [1,N]}$ is an orthonormal basis of $\R^N$. We set $C=\Pi_{i=1}^N (0,L_i)$ the periodicity cell. Fixing $(x_{1},...,x_{k-1},x_{k+1},...,x_{N})$, one can rearrange the function $x_{k}\mapsto \mu (x_{1},...,x_{k},...,x_{N})$. This is called the \textit{Steiner symmetrization} of the function $\mu$ in $x_{k}$. Performing successive rearrangements with respect to $x_{1},...,x_{N}$, one obtains a periodic function $\mu^{*}$ which is symmetric with respect to the planes $x_{k}=L_{k}/2$, nondecreasing in $x_{k}$ on the set $\{x_{k}\in (0,L_{k}/2)\}$, with the same distribution function as $\mu$. We underline that these conditions do not give a unique function $\mu^{*}$ and the way the symmetrization is carried out can lead to different symmetrized functions. In the sequel, we will call $\mu^{*}$ the function that is obtained after successive rearrangements in $x_{1},...,x_{N}$.

Using the Rayleigh characterization (\ref{RayleighSteiner}) and the properties of the Steiner symmetrization, one can prove (see \cite{Base1}) that $k_{0}(\mu^{*})\leq k_{0}(\mu)$. But this is not true anymore if one considers the periodic principal eigenvalue of the nonsymmetric operators $L_{\lambda e}$ defined by (\ref{LlambdamultiD}):

\begin{prop} \label{cexSteiner}
 If $N\geq 2$, then there exist $e\in\mathbb{S}^{N-1}$, $\lambda>0$ and $\mu$ such that 
\[k_{\lambda e}(\mu^*)>k_{\lambda e}(\mu) \hbox{ and } c^*_e(\mu^*)<c^*_e(\mu).\]
\end{prop}

In dimension $N\geq 2$, there exists many functions $\mu$ which have the same distribution function as $\mu$ and which are Steiner symmetric. This seems to be the main reason why Theorem \ref{influencesteiner} does not hold true. 

But we know that a minimizer exists: it has been proved by the author in \cite{eigenvalue} that for all $\lambda>0$, given $\mu^+>0$ and $m\in (0,1)$, the following infimum
$$\inf_{|\Omega|=m|C|}k_{\lambda e}(\mu^+ 1_\Omega)$$ 
is achieved for some measurable set $\Omega_{opt}\subset \overline{C}$. Moreover, if we define 
\[F = \big\{ \mu \in \mathcal{C}^0_{per}(\R^N), \ 
0 \leq \mu(x) \leq \mu^+ \ \hbox{for all} \ x, \frac{1}{|C|}
\int_{C} \mu = m\mu^+ \big\},\]
then the following infimums are equal
$$\inf_{\mu\in F}k_{\lambda e}(\mu)=\inf_{|\Omega|=m|C|}k_{\lambda e}(\mu^+ 1_\Omega).$$
In other words the ``bang-bang'' zero-order terms, that is the functions $\mu=\mu^+ 1_\Omega$, are minimizers in a wider class of functions. Studying the minimization problem for these functions is thus quite relevant. 

Proposition \ref{cexSteiner} yields that the minimizing set $\Omega_{opt}$ is not necessarily associated with the Steiner symmetrization of a given set $\Omega$, but this minimizing set exists and it would be relevant to try to determine its shape. Even in the case $\lambda=0$, we know that this set is symmetric and connected using the Steiner symmetrization, but L. Roques and F. Hamel \cite{optimal} have proved that it depends on the size $m$ and on $\mu^+$. In high dimension, more symmetries occur than in $\R$ and the Steiner symmetrization only takes into account the symmetries with respect to the planes $x_k=0$, $k=1,...,N$, but not the rotationnal symmetry for example. The numerical simulations of L. Roques and F. Hamel seem to show that $\Omega_{opt}$ is either a ball, a stripe or the complementary of a ball. But it is still an open problem to prove this conjecture. 

Let us also mention the works of L. Roques and R. Stoica \cite{RoquesStoica}. They used a stochastic modeling of fragmentation for a discrete approximation 
of our problem when $\mu$ has the form (\ref{patchy}). Their numerical computations proved that the discrete approximation of the periodic principal eigenvalue is a decreasing function of a given fragmentation parameter. These numerical computations show that there might exist a good notion of ``fragmentation'' in dimension $2$, which could make possible a generalization of Theorem \ref{influencesteiner}. 

Lastly, the counterexample given by Proposition \ref{cexSteiner} holds when $e$ is not a vector of the orthonormal basis. Then it strongly uses the fact that the successive Schwarz rearrangements are performed in the directions of the orthonormal basis. The case $e=e_1$ for example seems to be very different and in this case Theorem \ref{influencesteiner} may hold. We leave this problem as an open question.


\section{A new characterization of the periodic principal eigenvalue of a nonsymmetric operator}

\subsection{Statement of the characterization}

Theorem \ref{influencesteiner} is a corollary of a new formula for the periodic principal eigenvalue. This formula is available for general elliptic operators in any dimension.
We thus consider in this section the periodic principal eigenvalue $k_{\lambda e}(A,q,\mu)$ associated with the operator defined for all $\phi\in\mathcal{C}^2(\R^N)$ by
\begin{equation}-L_{\lambda}\phi=- \nabla \cdot(A(x) \nabla\phi) - 2\lambda eA(x)\nabla\phi -q(x)\cdot\nabla\phi-(\lambda^2 eA(x)e +\lambda\nabla\cdot(A(x)e)+\lambda q(x)\cdot e+\mu(x))\phi,\end{equation}
where $A$ is continuous, $q$ is of class $\mathcal{C}^1$ and these functions are both periodic in $x$:
\begin{equation} 
A(x+L_{i}\epsilon_i)=A(x) \hbox{ and } q(x+L_{i}\epsilon_i)=q(x) \hbox{ for all } i\in [1,N].
\end{equation}
Moreover, we assume that $A$ is a uniformly elliptic and continuous matrix field. That is,  
there exist some positive constants $\gamma$ and
$\Gamma$ such that for all $\xi \in \mathbb{R}^{N}, (x,t) \in
\mathbb{R}^{N} \times \mathbb{R}$ one has:
\begin{equation} \label{ellipticity} \gamma \|\xi\|^{2} \leq \sum_{1 \leq i,j \leq N} a_{i,j}(x,t) \xi_{i} \xi_{j} \leq \Gamma \|\xi\|^{2} \end{equation}
where $\|\xi\|=(|\xi_{1}|^{2}+...+|\xi_{N}|^{2})^{1/2}$.
Under these hypotheses, we know (see \cite{Pinsky}) that $k_{\lambda e}(A,q,\mu)$ is uniquely defined by the existence of a function $\phi$ such that 
\begin{equation} \label{eigen0bis}\left\{ \begin{array}{l}
 L_\lambda\phi = k_{\lambda e}(A,q,\mu)\phi,\\
 \phi > 0 \ \hbox{in} \ \mathbb{R}^N,\\
 \phi \hbox{ is periodic}.\\
\end{array} \right. \end{equation}

\bigskip

Our first formula holds when $q\equiv 0$. It involves the \textit{effective diffusivity} of a diffusion matrix $A$ in a given direction $e\in\mathbb{S}^{N-1}$, which is defined (see \cite{BensoussanLions}) by 
\begin{equation} \label{effectivediffusivity}  D_e(A)=\min_{\chi \in\mathcal{C}^1_{per}(\R^N)} \frac{1}{|C|}\int_{C} (e+\nabla\chi)A(x)(e+\nabla\chi). \end{equation} 

\begin{thm} \label{newformula}
One has
\begin{equation} \label{formuladim1} 
k_{\lambda e}(A,0,\mu)=\min_{\alpha\in \mathcal{A}} \Big(\int_C \nabla\alpha A(x)\nabla \alpha -\int_C \mu(x)\alpha^2 -\lambda^2 |C| D_e (\alpha^2 A) \Big),\end{equation}
where the effective diffusivity of a matrix field $D_e(\alpha^2 A)$ is defined by (\ref{effectivediffusivity}) and
\[\mathcal{A}=\{ \alpha\in \mathcal{C}^1_{per}(\R^N), \ \alpha>0, \ \int_C \alpha^2=1\}.\]
\end{thm}

This formula is deduced from the reformulation of a characterization of the principal eigenvalue of a nonsymmetric operator that has first been found by M. Donsker and S. Varadhan \cite{DonskerVaradhan} and by C. J. Holland \cite{Holland} in the case of Dirichlet or Neuman boundary conditions. The proofs used in these papers may certainly work in a periodic framework, but we give here an alternative proof and a reformulation that characterizes the periodic principal eigenvalue of a non-symmetric operator as the maximum of a family of periodic principal eigenvalues associated with self-adjoint operators:

\begin{thm} \label{Hollandthm}
For all $(A,q,\mu)$, one has
\begin{equation} \label{Hollandformula}
k_0(A,q,\mu)=\max_{\beta\in\mathcal{C}^1_{per}(\R^N)} k_0(A,0,\mu+\nabla\beta A \nabla\beta+q\cdot\nabla\beta-\frac{\nabla\cdot q}{2}).\end{equation}
\end{thm}

One can check that this formula is satisfied in some particular cases. 
If $q\equiv 0$, then $\max_{\beta\in\mathcal{C}^1_{per}(\R^N)} k_0(A,0,\mu+\nabla\beta A \nabla\beta)=k_0(A,0,\mu)$ and formula (\ref{Hollandformula}) is satisfied. More generally, if $q=A\nabla Q$, with $Q\in\mathcal{C}^1_{per}(\R^N)$, then one can check that this maximum is reached when $\beta=-Q/2$ and thus we find the classical formula 
\[k_0(A,A\nabla Q,\mu)=k_0(A,0,\mu-\frac{1}{4}\nabla Q A \nabla Q-\frac{1}{2}\nabla\cdot (A\nabla Q)).\]

The principal eigenvalues $k_0(A,0,q\cdot \nabla\beta + \nabla\beta A \nabla\beta +\mu -\frac{\nabla\cdot q}{2})$ are associated with self-adjoint operators. Thus they can be expressed using the Rayleigh characterization (\ref{RayleighSteiner}), which gives:
\begin{equation}
k_0(A,q,\mu)=\max_{\beta\in\mathcal{C}^1_{per}(\R^N)}\min_{\alpha\in\mathcal{A}}\big(\int_C \nabla\alpha A \nabla \alpha-\int_C \mu\alpha^2-\int_C (\nabla\beta A\nabla\beta+q\cdot \nabla\beta -\frac{\nabla\cdot q}{2})  \alpha^2\big).\end{equation}
By proving that it is possible to reverse the maximum and the minimum, we find the formula (\ref{infsup}), that was proved in \cite{DonskerVaradhan, Holland} for Neumann or Dirichlet boundary conditions:

\begin{prop} \label{Hollandprop}
For all $(A,q,\mu)$, one has:
\begin{equation} \label{infsup}\begin{array}{l} k_0(A,q,\mu)\\
=\min_{\alpha\in\mathcal{A}}\max_{\beta\in\mathcal{C}^1_{per}(\R^N)}\big(\int_C \nabla\alpha A \nabla \alpha-\int_C \mu\alpha^2-\int_C (\nabla\beta A\nabla\beta+q\cdot \nabla\beta -\frac{\nabla\cdot q}{2})  \alpha^2\big).\\
 \end{array}\end{equation}
\end{prop}

This last formula is the exact generalization of the formula proved by M. Donsker and S. Varadhan \cite{DonskerVaradhan} and by C. J. Holland \cite{Holland} in
bounded domains with Dirichlet or Neumann boudary conditions. This formula has been generalized to the principal eigenvalues of an elliptic operator 
with undefinite weight by T. Godoy, J.-P. Gossez and S. Paczka \cite{Godoyantimax} and to parabolic operators by T. Godoy, U. Kaufmann and S. Paczka \cite{Godoyparabolic}.

We now prove Theorem \ref{newformula}, Theorem \ref{Hollandthm} and Proposition \ref{Hollandprop}.


\subsection{Properties of the periodic principal eigenvalue}

We first remind some well-known properties of the periodic principal eigenvalues, that have been proved in \cite{Base1, Pinsky} for example.

\begin{prop} \label{classicalPinsky} \cite{Pinsky}
For all $A$, define 
\[G: \mu\in\mathcal{C}^0_{per}(\R^N)\mapsto k_0(A,0,\mu).\]
Then the following statements hold:

1) $G$ is concave,

2) $G$ is analytic and for all $\eta$, $G'_\mu(\eta)=-\int_C \eta \phi^2$,

where $\phi\in\mathcal{C}^2_{per}(\R^N)$ is the unique positive solution of
\[-\nabla\cdot (A(x)\nabla \phi)-\mu(x)\phi=k_0(A,0,\mu)\phi\]
normalized by $\int_C\phi^2=1$.
\end{prop}

We now define for all $(A,q,\mu)$:
\[F: \beta\in\mathcal{C}^1_{per}(\R^N)\mapsto k_0(A,0,\mu+\nabla\beta A \nabla\beta+q\cdot\nabla\beta-\frac{\nabla\cdot q}{2}).\]
In order to prove Theorem \ref{Hollandthm}, we have to study the maximization problem for the function $F$. 
We immediately get from Proposition \ref{classicalPinsky} the following result:

\begin{prop} \label{csqPinsky}
For all $(A,q,\mu)$, the following statements hold:

1) $F$ is concave,

2) $F$ is analytic and for all $h\in\mathcal{C}^1_{per}(\R^N)$, $F'_\beta(h)=-\int_C (2\nabla\beta A\nabla h+q\cdot\nabla h) \alpha^2$,

where $\alpha$ is the unique solution of 
\[-\nabla\cdot (A(x)\nabla\alpha) - (\mu-\frac{\nabla\cdot q}{2}+\nabla\beta A(x)\nabla\beta+q(x)\cdot \nabla\beta )\alpha=F(\beta)\alpha,\]
normalized by $\int_C\alpha^2=1$.
\end{prop}


\subsection{A self-adjoint equivalent problem}

We now prove that, up to some change of variable, the problem of finding a periodic principal eigenfunction of a nonsymmetric operator is equivalent to the problem of finding a periodic principal eigenfunction of a well-chosen self-adjoint operator. This change of variable has been introduced by S. Heinze in \cite{Heinze2005} and involves the periodic principal eigenfunction associated with the adjoint operator, which is uniquely defined by:
\begin{equation} \label{eigenadjoint0}\left\{
\begin{array}{l}
-\nabla \cdot (A(x)\nabla \tilde{\phi})+q(x)\cdot \nabla \tilde{\phi} -(\mu(x)-\nabla\cdot q(x))\tilde{\phi}=k_0(A,q,\mu)\tilde{\phi},\\
\tilde{\phi}>0, \ \tilde{\phi} \ \hbox{is periodic}, \ \int_C\phi\tilde{\phi}=1.\\
\end{array}\right. \end{equation}

\begin{prop} \label{propcdv}
Doing the change of variables
\begin{equation} \label{cdv} \alpha = \sqrt{\phi \tilde{\phi}}, \ \beta = \frac{1}{2} \ln (\frac{\phi}{\tilde{\phi}})\end{equation}
yields that finding solutions $(\phi,\tilde{\phi})\in\mathcal{C}^2(\R^N)\times\mathcal{C}^2(\R^N)$ of (\ref{eigen0bis})-(\ref{eigenadjoint0}) (with $\lambda=0$) is equivalent to finding a solution $(\alpha,\beta)\in\mathcal{C}^2(\R^N)\times\mathcal{C}^2(\R^N)$ of
\begin{equation} \label{pbmcdv}\left\{
\begin{array}{l}
-\nabla\cdot (A(x)\nabla\alpha) - (\mu-\frac{\nabla\cdot q}{2}+\nabla\beta A(x)\nabla\beta+q(x)\cdot \nabla\beta)\alpha=k_0(A,q,\mu)\alpha,\\
-\nabla\cdot \big( \alpha^2 (A(x)\nabla\beta+\frac{q(x)}{2})\big)=0,\\
\alpha>0, \ \alpha \ \hbox{is periodic},\\
\beta \ \hbox{is periodic}.\\
\end{array}\right.\end{equation}
\end{prop}

{\bf Remark:} The inverse change of variables is $\phi=\alpha e^\beta$ and $\tilde{\phi}=\alpha e^{-\beta}$. 

\bigskip

\noindent {\bf Proof.} We compute:
\begin{equation} \begin{array}{rcl}
\nabla\beta&=& \frac{1}{2} \big(\disp\frac{\nabla\phi}{\phi}-\disp\frac{\nabla\tilde{\phi}}{\tilde{\phi}}\big)\\
&&\\
\nabla\alpha &=& \disp\frac{\phi\nabla\tilde{\phi}}{2\alpha}+\disp\frac{\tilde{\phi}\nabla\phi}{2\alpha}\\
&&\\
-\nabla\cdot (A\nabla\alpha)&=& -\disp\frac{\phi}{2\alpha}\nabla\cdot (A\nabla\tilde{\phi})-\disp\frac{\tilde{\phi}}{2\alpha}\nabla\cdot(A\nabla\phi)-\disp\frac{\nabla\phi A\nabla\tilde{\phi}}{\alpha}\\
&& \\
&&+\disp\frac{1}{4\alpha^3} (\phi\nabla\tilde{\phi}+\tilde{\phi}\nabla\phi)A(\phi\nabla\tilde{\phi}+\tilde{\phi}\nabla\phi)\\
&&\\
&=& -\disp\frac{\phi}{2\alpha}\nabla\cdot (A\nabla\tilde{\phi})-\disp\frac{\tilde{\phi}}{2\alpha}\nabla\cdot(A\nabla\phi)+\disp\frac{1}{4\alpha^3} (\phi\nabla\tilde{\phi}-\tilde{\phi}\nabla\phi)A(\phi\nabla\tilde{\phi}-\tilde{\phi}\nabla\phi)\\
&&\\
&=& -\disp\frac{\phi}{2\alpha}\nabla\cdot (A\nabla\tilde{\phi})-\disp\frac{\tilde{\phi}}{2\alpha}\nabla\cdot(A\nabla\phi)+\nabla\beta A\nabla\beta\alpha.\\
\end{array} \end{equation}
This gives
\begin{equation} \begin{array}{l}
-\nabla\cdot (A(x)\nabla\alpha) - (\mu-\frac{\nabla\cdot q}{2}+\nabla\beta A(x)\nabla\beta+q\cdot \nabla\beta)\alpha\\
\\
=-\disp\frac{\phi}{2\alpha}\nabla\cdot (A\nabla\tilde{\phi})-\disp\frac{\tilde{\phi}}{2\alpha}\nabla\cdot(A\nabla\phi)-(\mu-\frac{\nabla\cdot q}{2}+q\cdot \nabla\beta)\alpha\\
\\
=\disp\frac{\phi}{2\alpha}(q\cdot \nabla \tilde{\phi}+\mu\tilde{\phi}+k_0(A,q,\mu)\tilde{\phi})+\disp\frac{\tilde{\phi}}{2\alpha}(-q\cdot \nabla \phi+\mu\phi-\nabla\cdot q \phi+k_0(A,q,\mu)\phi)\\
- (\mu-\frac{\nabla\cdot q}{2}+q\cdot \nabla\beta)\alpha\\
\\
=\disp\frac{\phi}{2\alpha}q\cdot \nabla \tilde{\phi}-\disp\frac{\tilde{\phi}}{2\alpha}q\cdot \nabla \phi+k_0(A,q,\mu)\alpha-q\cdot \nabla\beta\alpha\\
\\
=k_0(A,q,\mu)\alpha.\\
\end{array} \end{equation}
One can also compute:
\begin{equation} \begin{array}{l}
-\nabla\cdot \big( \alpha^2 (A(x)2\nabla\beta+q(x))\big)=-\nabla\cdot \big( \phi\tilde{\phi} (A(x)\disp\frac{\nabla\phi}{\phi}-A(x)\disp\frac{\nabla\tilde{\phi}}{\tilde{\phi}}+q(x))\big)\\
\\
=-\nabla\cdot \big( A(x)\tilde{\phi}\nabla\phi-A(x)\phi\nabla\tilde{\phi}+q(x)\tilde{\phi}\phi\big)\\
\\
=q(x)\tilde{\phi}\nabla\phi+(\mu+k_0(A,q,\mu))\tilde{\phi}\phi+q(x)\phi\nabla\tilde{\phi}-(\mu-\nabla\cdot q+k_0(A,q,\mu))\tilde{\phi}\phi\\
-(\nabla\cdot q(x)) \tilde{\phi}\phi-q(x)\nabla (\tilde{\phi}\phi)\\
\\
=q(x)\tilde{\phi}\nabla\phi+q(x)\phi\nabla\tilde{\phi}-q(x)\nabla (\tilde{\phi}\phi)=0.\\
\end{array} \end{equation}

In the other hand, if $\alpha$ and $\beta$ are given solutions of (\ref{pbmcdv}), then setting $\phi=\alpha e^\beta$ and $\tilde{\phi}=\alpha e^{-\beta}$, 
the previous equalities yield that $\phi$ and $\tilde{\phi}$ satisfy (\ref{eigen0bis}) with $\lambda=0$ and (\ref{eigenadjoint0}).$\Box$


\subsection{Proof of Theorem \ref{Hollandthm} and Proposition \ref{Hollandprop}}

Gathering Propositions \ref{csqPinsky} and \ref{propcdv}, we are able to conclude the proof of Theorem \ref{Hollandthm}.

\bigskip

\noindent {\bf Proof of Theorem \ref{Hollandthm}.} Consider $(\phi,\tilde{\phi})\in\mathcal{C}^2(\R^N)\times\mathcal{C}^2(\R^N)$ the solutions of (\ref{eigen0bis}-\ref{eigenadjoint0}) with $\lambda=0$ and $(\alpha,\beta)$ defined by (\ref{cdv}). Then Propositions \ref{csqPinsky} and \ref{propcdv} yield that for all $h\in\mathcal{C}^1_{per}(\R^N)$:
\[F'_{\beta}(h)=-\int_C (2\nabla\beta A\nabla h+q\cdot\nabla h) \alpha^2=\int_C  \nabla\cdot\big(\alpha^2(2A\nabla\beta+q)\big) h=0.\]
Thus $\beta$ is a critical point of $F$. As $F$ is concave from Proposition \ref{csqPinsky}, the function $F$ reaches its maximum in $\beta$. 

One also knows from Proposition \ref{propcdv} that $\alpha$ is the periodic principal eigenfunction associated with $F(\beta)$. This gives $F(\beta)=k_0(A,q,\mu)$. Thus Theorem \ref{Hollandthm} is proved. $\Box$

\bigskip

\noindent {\bf Proof of Proposition \ref{Hollandprop}.} Define 
\[J(\alpha',\beta')=\int_C \nabla\alpha' A(x) \nabla \alpha'-\int_C (\mu(x)+\nabla\beta' A(x)\nabla\beta'+q(x)\cdot\nabla \beta'-\frac{\nabla\cdot q}{2})  \alpha'^2.\]
We know from Theorem \ref{Hollandthm} and from the Rayleigh characterization of the periodic principal eigenvalue of a self-adjoint operator that 
\[k_0(A,q,\mu)=\max_{\beta'\in\mathcal{C}^1_{per}(\R^N)}F(\beta')=\max_{\beta'\in\mathcal{C}^1_{per}(\R^N)}\min_{\alpha'\in\A} J(\alpha',\beta').\]
Thus, in order to prove (\ref{infsup}), it is only left to prove that 
\begin{equation} \label{saddlenodeJ} \max_{\beta'\in\mathcal{C}^1_{per}(\R^N)}\min_{\alpha'\in\A} J(\alpha',\beta')=\min_{\alpha'\in\A}\max_{\beta'\in\mathcal{C}^1_{per}(\R^N)} J(\alpha',\beta').\end{equation}
Consider $(\alpha,\beta)$ as in the proof of Theorem \ref{Hollandthm}. 
As $\alpha$ is the periodic principal eigenfunction associated with $F(\beta)=k_0(A,0,\mu-\frac{\nabla\cdot q}{2}+\nabla\beta A\nabla\beta+q\cdot\nabla\beta)$, the Rayleigh characterization of this principal eigenvalue gives:
\[\forall \alpha'\in\mathcal{A}, \ J(\alpha,\beta) \leq J(\alpha',\beta).\]
On the other hand, we know that 
\[\begin{array}{l} \max_{\beta'\in\mathcal{C}^1_{per}(\R^N)}J(\alpha,\beta')\\
\\
=\max_{\beta'\in\mathcal{C}^1_{per}(\R^N)}\int_C \nabla\alpha A(x) \nabla \alpha-\int_C (\mu(x)+\nabla\beta' A(x)\nabla\beta'+q(x)\cdot\nabla \beta'-\frac{\nabla\cdot q}{2}) \alpha^2 \\
\\
=\int_C \nabla\alpha A(x) \nabla \alpha-\min_{\beta'\in\mathcal{C}^1_{per}(\R^N)}\big(\int_C (\mu(x)+\nabla\beta' A(x)\nabla\beta'+q(x)\cdot\nabla \beta'-\frac{\nabla\cdot q}{2}) \alpha^2\big)\\
\\
= J(\alpha,\beta_0),\\ \end{array} \]
where, by differentiation of $\beta'\mapsto \int_C (\nabla\beta' A(x)\nabla\beta'+q(x)\cdot\nabla \beta')$, the function $\beta_0$ is the unique solution (up to the addition of a constant) of
\begin{equation} \label{saddlenode1}
\left\{\begin{array}{l}  -\nabla \cdot \big(\alpha^2 A(x) (\nabla\beta_{0}+\frac{q(x)}{2})\big)=0,\\
\beta_{0} \ \hbox{is periodic}.\\
\end{array}\right.  \end{equation}
But Proposition \ref{propcdv} yields that $\beta$ is a solution (\ref{saddlenode1}). Thus there exists a constant $C$ such that $\beta_0=\beta+C$ and then $J(\alpha,\beta)=\max_{\beta'\in\mathcal{C}^1_{per}(\R^N)}J(\alpha,\beta')$.

This finally gives that $(\alpha,\beta)$ is a saddle node of the function $J$ and thus (\ref{saddlenodeJ}) holds. This ends the proof. $\Box$


\subsection{Proof of Theorem \ref{newformula}}

We now use Proposition \ref{Hollandprop} to compute $k_{\lambda e}(A,0,\mu)$.

\noindent {\bf Proof of Theorem \ref{newformula}.} As $k_{\lambda e}(A,0,\mu)=k_0(A,2\lambda e A, \mu + \lambda^2 eAe +\lambda \nabla\cdot (Ae))$, Proposition \ref{Hollandprop} gives (replacing $\beta$ by $\lambda\beta$):
\[\begin{array}{rcl}
k_{\lambda e}(A,0,\mu)&=&\min_{\alpha\in\mathcal{A}} \max_{\beta} \int_C \big(\nabla\alpha A\nabla\alpha-(\mu +\lambda^2\nabla \beta A \nabla\beta +\lambda^2 eAe +2\lambda^2 e A\nabla\beta)\alpha^2\big)\\
\\
&=&\min_{\alpha\in\mathcal{A}} \big(\int_C  \nabla\alpha A\nabla\alpha-\int_C\mu\alpha^2- \lambda^2\min_{\beta}\int_C (\nabla \beta+e)\alpha^2 A (\nabla \beta+e)\big)\\
\\
&=&\min_{\alpha\in\mathcal{A}} \big(\int_C \nabla\alpha A\nabla\alpha-\int_C\mu\alpha^2- \lambda^2|C| D_e(\alpha^2 A)\big).\\
\end{array}\] 
$\Box$


\section{Proof of Theorem \ref{influencesteiner} and construction of the counterexample in higher dimension}

\noindent {\bf Proof of Theorem \ref{influencesteiner}.} We use Proposition \ref{Hollandprop}. Set $C=(0,L)$. It is well-known that, in dimension $1$, for all positive function $a$ and for $e=\pm 1$, $D_e(a)^{-1}=\frac{1}{L} \int_0^L \frac{dx}{a(x)}$ (see \cite{BensoussanLions} for example). 
Thus
\begin{equation} \label{formuladim1-2}
k_{\lambda}(\mu)=\min_{\alpha\in \mathcal{A}} \Big(\int_0^L \alpha'^2-\int_0^L \mu\alpha^2-\displaystyle{\frac{L\lambda^2}{\frac{1}{L}\int_0^L \frac{1}{\alpha^2}} }\ \Big).
\end{equation}

Consider now any $\alpha\in\A$ and its Steiner rearrangement $\alpha^*$. As $\alpha$ and $\alpha^*$ have the same distribution function, one has 
\[\int_0^L(\alpha^*)^2=\int_0^L\alpha^2=1 \ \hbox{and} \ \int_0^L\frac{1}{(\alpha^*)^2}=\int_0^L\frac{1}{\alpha^2}.\]
It is well known (see \cite{Kawohl} for example) that
\begin{equation} \label{influencesteiner-proof0}\begin{array}{cc} \int_0^L (\alpha^*)'^2\leq \int_0^L \alpha'^2 &\hbox{ (Poly\`a inequality)},\\ \int_0^L \mu\alpha^2\leq \int_0^L \mu^*(\alpha^*)^2 &\hbox{ (Hardy-Littlewood inequality)}.\\ \end{array}\end{equation}
Thus it immediately follows from (\ref{formuladim1-2}) that $k_{\lambda}(\mu^*)\leq k_{\lambda}(\mu)$. 
$\Box$

\bigskip

In order to prove Proposition \ref{cexSteiner}, we first prove a result of independant interest:
\begin{prop} \label{equivalenceconj}
For all $\mu_{1}, \mu_{2}\in\mathcal{C}_{per}^0(\R^N)$ and $e\in\mathbb{S}^{N-1}$, the following two assertions are equivalent:

1) for all $M>k_0(\mu_{1})$, one has: $c^{*}_e(\mu_{1}+M)\leq c^{*}_e(\mu_{2}+M)$,

2) for all $\lambda>0$, one has: $k_{\lambda e}(\mu_{1})\geq k_{\lambda e}(\mu_{2})$.
\end{prop}

Basically, this result proves that trying to show that $c^*_e(\mu_1)\leq c^*_e(\mu_2)$ by testing if $k_{\lambda e}(\mu_1)\geq k_{\lambda e} (\mu_2)$ is quite relevant.

\bigskip

\noindent {\bf Proof of Proposition \ref{equivalenceconj}.} If $2)$ is true, then $1)$ immediately follows from (\ref{formulec^*}). Assume that $2)$ is false now: there exists some $\lambda_{0}>0$ such that $k_{\lambda_{0} e}(\mu_{1})< k_{\lambda_{0} e}(\mu_{2})$. Set $c_{0}=-\partial_{\lambda}k_{\lambda e}(\mu_{1})\big|_{\lambda=\lambda_0}$. The function $\lambda\mapsto k_{\lambda e}(\mu_{1})+\lambda c_{0}$ reaches its maximum for $\lambda=\lambda_{0}$ since it is concave and its derivative vanishes for $\lambda=\lambda_{0}$. Set:
\[M=k_{\lambda_{0} e}(\mu_{1})+\lambda_{0} c_{0}>k_{0}(\mu_{1}).\]
One has $c^{*}(\mu_{1}+M)=c_{0}$ since 
\[\max_{\lambda>0}(k_{\lambda e}(\mu_{1}+M)+\lambda c_0)=k_{\lambda_{0} e}(\mu_{1}+M)+\lambda_{0} c_0=0.\]
In the other hand:
\[k_{\lambda_{0} e}(\mu_{2}+M)+\lambda_{0}c_{0}=k_{\lambda_{0} e}(\mu_{2})-k_{\lambda_{0} e}(\mu_{1})>0,\]
thus:
\[c^{*}(\mu_{1}+M)>c^{*}(\mu_{2}+M),\]
which contradicts $1)$. The equivalence is proved. $\Box$

\bigskip

\noindent {\bf Proof of Proposition \ref{cexSteiner}.} We set $j_{\lambda e}(\mu)=k_{\lambda e}(\mu)+\lambda^{2}$ for all $\lambda>0$. It has been proved in \cite{Largedrift} that:
\begin{equation} \label{cvlargedrift} j_{\lambda e}(\mu)\rightarrow \min_{e\cdot\nabla\phi=0}\frac{\int_{C}(|\nabla\phi|^{2}-\mu\phi^{2})}{\int_{C}\phi^{2}} \hbox{ as } \lambda\rightarrow+\infty. \end{equation}

For $N=2$, take $e = (1,-1)$ and $\mu(x,y)=\mu_{0}(x+y)$, where
$\mu_{0}^{*}=\mu_{0}$ and $\mu_{0}$ is not constant. Take
$\varphi_{0}$ defined in dimension $1$ by:
\[ \left\{ \begin{array}{l}
-\varphi_{0}''- \frac{\mu_{0}}{2} \varphi_{0} =
k_{0}({\frac{\mu_{0}}{2}}) \varphi_{0},\\
\varphi_{0} > 0,\\
\varphi_{0} \ \hbox{L-periodic}.\\ \end{array} \right. \]

Set $\varphi(x,y)=\varphi_{0}(x+y)$, this function satisfies:
\[ (-\Delta \varphi - 2 \lambda e.\nabla \varphi - \mu \varphi)(x,y)
=(-2 \varphi_{0}'' - \mu_{0} \varphi_{0})(x+y)= 2
k_0(\frac{\mu_{0}}{2}) \varphi(x,y)\]
This gives $j_{\lambda e}(\mu)= 2 k_{0}(\frac{\mu_{0}}{2})$ which does not depend on $\lambda$.
Rearranging first in $x$ and then in $y$, one gets $\mu^{*}(x,y)=\mu_{0}(x)$. Using $(\ref{cvlargedrift})$, one gets:
\[j_{\lambda e}(\mu^{*}) \rightarrow \min_{\partial_{x} \psi =\partial_{y} \psi}
\frac{\int_{C}2|\partial_{x} \psi|^{2}- \mu_{0}(x)
\psi^{2}}{\int_{C}\psi^{2}} \hbox{ as } \lambda\rightarrow+\infty.\] 
Take $\psi$ a function for which this minimum is reached.
Then for all $y \in [0,L]$, setting $\psi_{y}(x)=\psi(x,y)$, one gets:
\begin{equation} \label{cex-1} 2 \int_{0}^{L} (|\psi_{y}'|^{2}-\frac{\mu_{0}(x)}{2}
\psi_{y}^{2}) \geq 2 k_0(\frac{\mu_{0}}{2}) \times
\int_{0}^{L} \psi_{y}^{2}\end{equation} 
using the definition of
$k_0(\frac{\mu_{0}}{2})$. Integrating in $y\in (0,L)$, this yields that:
\[\lim_{\lambda\rightarrow+\infty} j_{\lambda e}(\mu^{*}) \geq 2 k_0(\frac{\mu_{0}}{2}) = \lim_{\lambda\rightarrow+\infty} j_{\lambda e}(\mu)\]
If the equality holds, then for all $y$ (\ref{cex-1}) is an equality. But then $\psi_{y}$ is some eigenfunction associated with $k_{0}(\frac{\mu_{0}}{2})$ and thus it can be written $\psi_{y}(x)=\psi(x,y)=\theta(y)\varphi_{0}(x)$ for all $x,y$. But $\partial_{y}\psi=\partial_{x}\psi$, which gives that $\theta$ and $\varphi_{0}$ are constant. This is a contradiction since $\mu_{0}$ is not constant. 

Thus $\lim_{\lambda\rightarrow+\infty} j_{\lambda e}(\mu^{*})>\lim_{\lambda\rightarrow+\infty} j_{\lambda e}(\mu)$ and 
$$\hbox{there exists } \lambda>0 \hbox{ such that } k_{\lambda e}(\mu^*)>k_{\lambda e}(\mu).$$
If $N>2$, setting $\mu(x) = \mu_{0}(x_{1}+ x_{2})$ and $e = (1,-1,0,...,0)$, we are back to the case $N = 2$. $\Box$


\section{Influence of the period}\label{sectionperiod}

One can wonder if an analogue of Theorem \ref{influencesteiner} still holds with heterogeneous diffusion and advection coefficients. 
Dependence relation between the periodic principal eigenvalue and the shape of the diffusion coefficients is not well-known. 

Thus, we quantify in this section the ``fragmentation'' of the environment by investigating the dependence between $k_{\lambda e}(A,q,\mu)$ and the period of the coefficients. 
Define 
\[A_L(x)=A(\frac{x}{L}), \ q_L(x)=q(\frac{x}{L}), \ \mu_L(x)=\mu(\frac{x}{L}).\]
If $(A,q,\mu)$ is $1$-periodic, then $(A_L,q_L,\mu_L)$ is $L$-periodic. In \cite{Shigesada2}, N. Shigesada, K. Kawasaki and E. Teramoto carried out some numerical computations which seem to show that $L\mapsto c^*_e(A_L,0,\mu_L)$ is a nondecreasing function. We give here a rigourous proof of this fact and asymptotics when $L\rightarrow 0$:
 
\begin{prop} \label{influenceperiod}
For all $(A,q,\mu)$ with $\nabla\cdot q=0$, the function $L\mapsto k_{\lambda e}(A_L,q_L,\mu_L)$ is a nonincreasing function. Moreover, if $\int_{C}\mu> 0$ and $\nabla\cdot q\equiv 0$, one has 
\[ k_{\lambda e}(A_L, q_L,\mu_L)\rightarrow -\overline{\mu}-\lambda^2 D_e(A) \hbox{ as } L\rightarrow 0 \hbox{ and}\]
\[c^*_e(A_L,q_L,\mu_L)\rightarrow 2\sqrt{ D_e(A) \overline{\mu}} \hbox{ as } L\rightarrow 0.\]
where $\overline{\mu}=\frac{1}{|C|} \int_{C} \mu$.
\end{prop}

This immediately gives that $L\mapsto c^*_e(A_L,q_L,\mu_L)$ is a nondecreasing function. 
As an immediate corollary of Proposition \ref{influenceperiod}, one gets 
\[c^*_e(A_L,q_L,\mu_L)\geq 2\sqrt{ D_e(A) \overline{\mu}} \hbox{ for all } L>0.\]
This inequality has never been proved before as far as we know.

The asymptotics when $L\rightarrow 0$ have already been proved in some particular cases. 
For example, in dimension $1$, it is easy to compute $D_e(a)^{-1}=\frac{1}{|C|}\int_C \frac{1}{a}$. 
In this particular case, the convergence result of Proposition \ref{influenceperiod} (but not the monotonicity result) has already been proved by M. El Smaily, 
F. Hamel and L. Roques in \cite{EHR}. Similarly, it has been proved in \cite{Base2} when $A=I_N$ and in \cite{ElSmaily} when $\nabla\cdot (Ae)=0$. 
In these two cases, $D_e(A)=\frac{1}{|C|}\int_C eA(x)edx$. 

It is much more difficult to get asymptotics for $L\rightarrow +\infty$. F. Hamel, L. Roques and J. Favard have computed this limit when $N=1$ and $A$ and $\mu$ only take two values. 
We will generalize, together with F. Hamel and L. Roques, this convergence result to general $A$ and $\mu$ in a forthcoming work.

\bigskip

\noindent {\bf Proof if Proposition \ref{influenceperiod}.} A simple rescaling argument (see \cite{dependance} for example) gives for all $\lambda>0$:
\[k_{\lambda e}(A_L,q_L,\mu_L)=\frac{1}{L^2} k_{\lambda Le}(A,L q, L^2 \mu).\]
Using Theorem \ref{Hollandthm}, this gives
\begin{equation} \label{hom7} \begin{array}{l} 
k_{\lambda e}(A_L,q_L,\mu_L)\\
=\max_{\beta\in\mathcal{C}^1_{per}(\R^N)}\frac{1}{L^2} k_{0}(A,0, \nabla\beta A \nabla \beta +L q\cdot \nabla \beta +2L\lambda e A \nabla\beta+L^2\lambda e\cdot q +\lambda^2 L^2 eAe +  L^2 \mu)\\
=\max_{\beta_0\in\mathcal{C}^1_{per}(\R^N)}\frac{1}{L^2} k_{0}\Big(A,0, L^2\big(\lambda^2(\nabla\beta_0+e) A (\nabla \beta_0+e) +\lambda q\cdot (\nabla \beta_0+e)+\mu\big)\Big),\\
\end{array} \end{equation}
with the change of variable $\beta= L\lambda\beta_0$. As $\eta\mapsto k_0(A,0,\eta)$ is concave, this gives for all $L>1$:
\[ \begin{array}{l} 
k_{\lambda e}(A_L,q_L,\mu_L)\\
\leq  \max_{\beta_0\in\mathcal{C}^1_{per}(\R^N)} k_{0}\big(A,0, \lambda^2(\nabla\beta_0+e) A (\nabla \beta_0+e) +\lambda q\cdot (\nabla \beta_0+e)+\mu\big)= k_{\lambda e}(A,q,\mu).\\
\end{array} \]
This concludes the first part of the proof since $L>1$ and $(A,q,\mu)$ are general.

\bigskip

When $L\rightarrow 0$, (\ref{hom7}) gives 
\begin{equation} \label{hom0} \begin{array}{l} 
\lim_{L\rightarrow 0} k_{\lambda e}(A_L,q_L,\mu_L)\\
\geq \max_{\beta}\lim_{L\rightarrow 0} \frac{1}{L^2} k_{0}\Big(A,0, L^2\big(\lambda^2(\nabla\beta+e) A (\nabla \beta+e) +\lambda q\cdot (\nabla \beta+e)+\mu\big)\Big).\\
\end{array} \end{equation}
Moreover, for all $\beta$, using 2) of Proposition \ref{classicalPinsky} :
\begin{equation}\label{hom2}\begin{array}{l}
\lim_{L\rightarrow 0}  \frac{1}{L^2} k_{0}\Big(A,0, L^2\big(\lambda^2(\nabla\beta+e) A (\nabla \beta+e) +\lambda q\cdot (\nabla \beta+e)+\mu\big)\Big)\\
\\
= -\frac{1}{|C|}\int_C \big(\lambda^2(\nabla\beta+e) A (\nabla \beta+e) +\lambda q\cdot (\nabla \beta+e)+\mu\big)\\
\\
= -\frac{\lambda^2}{|C|}\int_C (\nabla\beta+e) A (\nabla \beta+e) -\overline{\mu},\\
\end{array} \end{equation}
since $\nabla\cdot q = 0$ and $\int_C q\cdot e=0$. Thus 
\begin{equation} \label{hom6} 
\lim_{L\rightarrow 0} k_{\lambda e}(A_L,q_L,\mu_L) \geq \max_\beta\big( -\frac{\lambda^2}{|C|}\int_C (\nabla\beta+e) A (\nabla \beta+e) -\overline{\mu}\big) = -\lambda^2 D_e(A) -\overline{\mu}.
 \end{equation}
On the other hand, as $k_0(A,0,\eta)\leq -\frac{1}{|C|}\int_C \eta$ for all $\eta$, one has for all $L$:
\begin{equation} \label{hom1} k_{\lambda e}(A_L,q_L,\mu_L)\leq -\min_{\beta\in\mathcal{C}^1_{per}(\R^N)}\frac{\lambda^2}{|C|}\int_C (\nabla\beta+e) A (\nabla \beta+e) -\overline{\mu}=-\lambda^2 D_e(A)-\overline{\mu}.\end{equation}

Gathering (\ref{hom6}) and (\ref{hom1}), one finally gets
\begin{equation} \lim_{L\rightarrow 0} k_{\lambda e}(A_L,q_L,\mu_L)=-\overline{\mu}-\lambda^2 D_e(A).\end{equation}
Moreover, this limit is locally uniform in $\lambda$ since $\lambda\mapsto k_{\lambda e}(A_L,q_L,\mu_L)$ is concave for all $L$.

\bigskip

Now, (\ref{hom1}) gives:
\begin{equation} \label{lowerboundc^*} c^*_e(A_L,q_L,\mu_L)=\min_{\lambda>0} \frac{-k_{\lambda e}(A_L,q_L,\mu_L)}{\lambda} \geq \min_{\lambda>0}\big(\frac{1}{\lambda} \overline{\mu}+\lambda D_e(A)\big)=2\sqrt{ \overline{\mu} D_e(A)}.\end{equation}

For all $L>0$, define $\lambda_L>0$ such that
\begin{equation} \label{deflambdaL} c^*_e(A_L,q_L,\mu_L)=\min_{\lambda>0}\frac{- k_{\lambda e}(A_L,q_L,\mu_L)}{\lambda}=\frac{- k_{\lambda_L e}(A_L,q_L,\mu_L)}{\lambda_L}.\end{equation}
We know from the maximum principle that
$$\begin{array}{rcl} k_{\lambda e}(A_L,q_L,\mu_L)&\geq &\min_C -(\lambda^2 eA_Le +\lambda\nabla\cdot (A_L e ) +\lambda q_L\cdot e+\mu_L)\\
&\geq& -\lambda^2\max_C(eA_Le)-\max_C(\mu_L)-\lambda\max_C |q_L|-\lambda\max_C  |\nabla\cdot (A_Le)|,\\ \end{array}$$
this gives
\[c^*_e(A_L,q_L,\mu_L)\leq 2\sqrt{ \max_C(eAe)\max_C(\mu)}+\max_C |q|+\max_C |\nabla\cdot (Ae)|\]
and 
\[k_{\lambda_L e}(A_L,q_L,\mu_L)\leq -\overline{\mu}-\lambda_L^2 D_e(A),\]
one has 
\[\lambda_L^2\leq \frac{1}{D_e(A)} (\overline{\mu}+ 2\sqrt{ \max_C(eAe)\max_C(\mu)}+\max_C |q|+\max_C |\nabla\cdot (Ae)|)\]
and thus the family $(\lambda_L)_{L>0}$ is bounded.
Hence, taking a sequence $L_n\rightarrow 0$, one may assume, up to extraction, that $\lambda_{L_n}\rightarrow \lambda_0\geq 0$. If $\lambda_0=0$, then we get from (\ref{deflambdaL}) 
\[c^*_e(A_L,q_L,\mu_L)\sim \frac{\overline{\mu}}{\lambda_L}\rightarrow-\infty,\]
which would contradict (\ref{lowerboundc^*}). Thus $\lambda_0>0$.
This gives
\begin{equation} \label{hom5}c^*_e(A_L,q_L,\mu_L)=\frac{- k_{\lambda_L e}(A_L,q_L,\mu_L)}{\lambda_L}\rightarrow \frac{1}{\lambda_0} \overline{\mu}+\lambda_0 D_e(A) \ \hbox{as} \ L\rightarrow 0.\end{equation}

We conclude from (\ref{lowerboundc^*}) and (\ref{hom5}) that
\begin{equation} \lim_{L\rightarrow 0} c^*_e(A_L,q_L,\mu_L)= \min_{\lambda> 0}\big( \frac{1}{\lambda} \overline{\mu}+\lambda D_e(A)\big)=2\sqrt{ \overline{\mu} D_e(A)}.\end{equation}
$\Box$


\section{Rearrangement of the diffusion coefficient}

One can wonder if an analogue of Theorem \ref{influencesteiner} still holds if one considers a heterogeneous diffusion coefficient. We thus consider the periodic principal eigenvalue $k_\lambda(a,\mu)$ associated with the modified operator defined for all $\phi$ by
\begin{equation}-\tilde{L_{\lambda}}\phi=- (a(x)\phi')' - 2\lambda a(x)\phi' -(\lambda^2 a(x) +\lambda a'(x)+\mu(x))\phi,\end{equation}
where $a$ is continuous, positive and $L-$periodic. We define $\displaystyle a_*= \frac{1}{(\frac{1}{a})^*}$.

\begin{thm} \label{diffusionthm}
Assume that $N=1$ and consider $(a,\mu)\in\mathcal{C}^0_{per}(\R)$ and $\lambda\in\R$ such that $a$ is positive. Then one has
\begin{equation} \label{dep-a} k_\lambda(a_*,\frac{(\mu a)^*}{a_*})\leq k_\lambda (a,\mu).\end{equation}
In particular, when $\mu=\mu_0>0$ is a constant, one has
\begin{equation} \label{dep-a2} k_\lambda(a_*,\mu_0)\leq k_\lambda (a,\mu_0) \hbox{ and } c^*(a_*,\mu_0)\geq c^*(a,\mu_0).\end{equation}
\end{thm}

This kind of inequality has never been proved before as far as we know, even in the self-adjoint case $\lambda=0$. If $a=1$, then we are back to Theorem \ref{influencesteiner}.

In the particular case of sinuisoidal diffusion and growth rate, when only the phase of the diffusion is varying, N. Kinezaki, K. Kawasaki and N. Shigesada \cite{Shigesadasinus} carried out some numerical computations which yield that the speed reaches its maximum when the diffusion and the growth have opposite phases. 
Trying to write this numerical results in a rigorous and general way, we can formulate the following conjecture:
\begin{equation} \label{conjecturediffusion} k_\lambda(a,\mu)\geq k_\lambda(a_*,\mu^*) \hbox{ for all }\lambda, a, \mu.\end{equation} 
The interpretation for this inequation is the following: we expect that the propagation speed reaches its maximum when there is a slow-down in the favourable area and a speed-up in the unfavourables ones. 

\bigskip

{\bf Proof of Theorem \ref{diffusionthm}.} First, if $(\ref{dep-a})$ is true for all $(a,\mu)$, then taking $\mu\equiv 0$ gives $k_\lambda(a_*,0)\leq k_\lambda (a,0)$. As $k_\lambda(a,\mu_0)=k_\lambda(a,0)-\mu_0$ for all constant $\mu_0$, we get (\ref{dep-a2}). 

Next, Proposition \ref{Hollandprop} yields
\[k_\lambda(a,\mu)= \min_{\psi>0} \displaystyle\frac{1}{\int_0^L \psi^2} \Big(\int_0^{L} a(x)\psi'^2dx - \int_0^L\mu(x)\psi^2dx - \frac{\lambda^2 L^2}{\int_0^L \frac{dx}{a(x)\psi^2}}\Big).\]
Consider some test-function $\psi>0$ and define $\phi=\psi \circ f_a^{-1}$, where $f_a(x)=\int_0^x \frac{dt}{a(t)}$. Let $L_1=f_a(L)$. As $a$ is positive, $f_a$ is an increasing bijection from $(0,L)$ to $(0,L_1)$. 
Thus, the change of variable $y=f_a(x)$ gives 
\begin{equation} \label{cdvdiffusion-1} k_\lambda(a,\mu)= \min_{\phi>0}\displaystyle\frac{1}{\int_0^{L_1} \rho(y)\phi^2}\Big(\int_0^{L_1} \phi'^2dy - \int_0^{L_1}\eta(y)\phi^2dy - \frac{\lambda^2 L_1^2}{\int_0^{L_1}\frac{dy}{\phi^2}}\Big),\end{equation}
where $\rho=a\circ f_a^{-1}$ and $\eta=(a\mu)\circ f_a^{-1}$. 

Taking $\phi^*$ as a test function and using the classical rearrangement inequalities, this gives
\begin{equation} \label{cdvdiffusion-2} k_\lambda(a,\mu)\geq \min_{\phi>0}\displaystyle\frac{1}{\int_0^{L_1} \rho_*(y)(\phi^*)^2}\Big(\int_0^{L_1} (\phi^*)'^2dy - \int_0^{L_1}\eta^*(y)(\phi^*)^2dy - \frac{\lambda^2 L_1^2}{\int_0^{L_1}\frac{dy}{(\phi^*)^2}}\Big).\end{equation}

We now prove that $\rho_*=a_*\circ f_{a_*}^{-1}$ and $\eta^*=(a\mu)^*\circ f_{a_*}^{-1}$.

\smallskip

{\bf A reverse change of variable.} We first construct some $a_0$ such that $\rho_*=a_0\circ f_{a_0}^{-1}$. 
Define $b$ the solution of
\begin{equation} \label{edo-cdv} b'(x)= \frac{1}{\rho_*(b(x))}, \ b(0)=0.\end{equation}
Set $a_0(x)=1/b'(x)$, so that $f_{a_0}(x)=b(x)$. This construction yields
\[\rho_*\circ f_{a_0}(x)=\rho_*(b(x))=\frac{1}{b'(x)}= a_0(x).\]
We now prove that $a_0\equiv a_*$ using the uniqueness of the Schwarz rearrangement. 

\smallskip

{\bf Equality of the distribution functions.}
For all $\lambda\in\R$, one has
\[mes \{a\geq\lambda\}= \int_0^L 1_{a(x)\geq \lambda}dx= \int_0^{L_1} \rho(y)1_{\rho(y)\geq\lambda}dy=\int_0^{L_1} \rho_*(y)1_{\rho_*(y)\geq\lambda}dy=mes \{a_0\geq\lambda\}\]
since $\rho$ and $\rho_*$ have the same distribution function. Thus $a_0$ and $a$ have the same distribution function.

\smallskip

{\bf Symmetry and periodicity.} The previous step gives $f_{a_0}(L)=\int_0^L\frac{1}{a_0}=\int_0^L\frac{1}{a}=L_1$. Thus $f_{a_0}$ and $x\mapsto L_1-f_{a_0}(L-x)$ are both solutions of (\ref{edo-cdv}) since $\rho_*$ is symmetric with respect to $L_1$. The Cauchy theorem yields that $f_{a_0}(x)=L_1-f_{a_0}(L-x)$ for all $x$. Deriving this equality, one gets that $a_0$ is symmetric with respect to $L/2$. This gives that $a_0$ is $L$-periodic.  

\smallskip

{\bf Monotonicity.}
As $f_{a_0}$ is increasing and $\rho_*$ is nonincreasing in $(0,L_1/2)$ and nondecreasing in $(-L_1/2,0)$, we know that the even function $a_0=\rho_*\circ f_{a_0}$ is nonincreasing in $(0,L_1/2)$ and nondecreasing in $(-L_1/2,0)$. 

\bigskip

As the Schwarz rearrangement of a function is uniquely defined, we conclude from all these properties that $a_0\equiv a_*$. Similarly, one can prove that $\eta^*=(a\mu)^*\circ f_{a_*}^{-1}$. Equation (\ref{cdvdiffusion-2}) gives
\begin{equation} \label{cdvdiffusion-3} k_\lambda(a,\mu)\geq \min_{\phi>0}\displaystyle\frac{1}{\int_0^{L_1} a_*\circ f_{a_*}^{-1}(y)(\phi^*)^2}\Big(\int_0^{L_1} (\phi^*)'^2dy - \int_0^{L_1}(a\mu)^*\circ f_{a_*}^{-1}(y)(\phi^*)^2dy - \frac{\lambda^2 L_1^2}{\int_0^{L_1}\frac{dy}{(\phi^*)^2}}\Big).\end{equation}
The change of variable $y=f_{a_*}(x)$ leads to
\begin{equation} \label{cdvdiffusion-4}\begin{array}{rcll} k_\lambda(a,\mu)&\geq& \min_{\phi>0} \displaystyle\frac{1}{\int_0^L (\phi^*\circ f_{a_*})^2} \Big(&\int_0^{L} a_*(x)(\phi^*\circ f_{a_*})'^2dx - \int_0^L(a\mu)^*(x)(\phi^*\circ f_{a_*})^2\frac{dx}{a_*(x)}\\
&&& - \lambda^2 L^2 \big(\int_0^L (a_*(x)(\phi^*\circ f_{a_*}))^{-2}dx\big)^{-1}\Big)\\
&&&\\
&\geq& k_\lambda(a_*,\frac{(a\mu)^*}{a_*}),&\\ \end{array} \end{equation}
using Proposition \ref{Hollandprop} and the test-function $\psi=\phi^*\circ f_{a_*}$. 
$\Box$

\bibliographystyle{plain}
\bibliography{biblio}

\end{document}